\documentclass{article}
\usepackage{amsmath,amssymb,latexsym}
\usepackage{graphicx}

\newtheorem{theorem}{Theorem}[section]

\newtheorem{corollary}[theorem]{Corollary}

\newtheorem{definition}[theorem]{Definition}

\newtheorem{lemma}[theorem]{Lemma}

\newtheorem{proposition}[theorem]{Proposition}

\newtheorem{remark}[theorem]{Remark}

\newenvironment{proof}[1][Proof]{\textbf{#1.} }{\ \rule{0.5em}{0.5em}}

\newcommand{\bth}{\begin{theorem}}
\newcommand{\eeth}{\end{theorem}}
\newcommand{\ble}{\begin{lemma}}
\newcommand{\ele}{\end{lemma}}
\newcommand{\bco}{\begin{corollary}}
\newcommand{\eco}{\end{corollary}}
\newcommand{\bde}{\begin{definition}}
\newcommand{\ede}{\end{definition}}
\newcommand{\bpr}{\begin{proposition}}
\newcommand{\epr}{\end{proposition}}
\newcommand{\bre}{\begin{remark}}
\newcommand{\ere}{\end{remark}}
\newcommand{\beq}{\begin{equation}}
\newcommand{\eeq}{\end{equation}}
\newcommand{\ben}{\begin{equation*}}
\newcommand{\een}{\end{equation*}}

\begin{document}
\title{\large\bf $W$-entropy formula for the Witten Laplacian on
manifolds with time dependent metrics and potentials}
\author{Songzi Li\thanks{Research partially supported by the China Scholarship Council.}, \ \ Xiang-Dong Li\thanks{Research supported by NSFC No. 11371351, Key Laboratory RCSDS, CAS, No. 2008DP173182, and a
Hundred Talents Project of AMSS,
CAS.}}

\maketitle

\begin{minipage}{120mm}
{\bf Abstract}. In this paper, we develop a new approach to prove the $W$-entropy formula for the Witten Laplacian via warped product on Riemannian manifolds and give a natural geometric interpretation of a quantity appeared in the $W$-entropy formula. Then we prove the $W$-entropy formula for the Witten Laplacian on compact Riemannian manifolds with time dependent metrics and potentials, and derive the $W$-entropy formula for the backward heat equation  associated with the Witten Laplacian on compact Riemannian manifolds equipped with Lott's modified Ricci flow. We also extend our results to complete Riemannian manifolds with negative $m$-dimensional Bakry-Emery Ricci curvature, and to compact Riemannian manifolds with $K$-super $m$-dimensional Bakry-Emery Ricci flow. As application, we prove that the optimal logarithmic Sobolev constant on compact manifolds equipped with the $K$-super $m$-dimensional Bakry-Emery Ricci flow is decreasing in time.
\end{minipage}


\section{Introduction}
Let $M$ be a complete Riemannian manifold with a fixed Riemnnian metric $g$ and a fixed potential $\phi\in C^2(M)$. Let  $d\mu=e^{-\phi}dv$, where $dv$ is the Riemannian volume measure on $(M, g)$. The Witten Laplacian, called also the weighted Laplacian,
\begin{eqnarray*}
L =\Delta -\nabla \phi\cdot\nabla \label{WL}
\end{eqnarray*}
is a self-adjoint and non-negative operator on $L^2(M, \mu)$. By It\^o's calculus, one can construct the symmetric diffusion process $X_t$ associated to the Witten Laplacian by solving the SDE
\begin{eqnarray*}
dX_t=\sqrt{2}dW_t-\nabla\phi(X_t)dt,
\end{eqnarray*}
where $W_t$ is the Brownian motion on $M$. Moreover, it is well known that the transition probability density function of the diffusion process $X_t$ is exactly the fundamental solution to the heat equation of $L$, i.e., the heat kernel of the Witten Laplacian $L$. In view of this, it is a fundamental problem to study the heat equation and the heat kernel of the Witten Laplacian on manifolds.

In recent years, important progress has been obtained in the study of the heat equation associated with the Witten Laplacian by using new ideas and new methods from geometric analysis, PDEs, and probability theory. In particular,  F. Otto \cite{Ot} introduced an infinite dimensional Riemannian structure on the Wasserstein space of probability measures on  $\mathbb{R}^n$ and proved that the heat equation
\begin{eqnarray}
\partial_t u=Lu\label{L-HE}
\end{eqnarray}
can be realized as the reverse gradient flow of the Bolztmann-Shannon entropy\footnote{Equivalently, the heat equation $(\ref{L-HE})$ is the gradient flow of $Ent(u)=-H(u)$ on the Wasserstein space $\mathcal{P}_2(\mathbb{R}^n)$ equipped with Otto's infinite dimensional Riemannian metric.}
$$
H(u)=-\int_M u\log ud\mu.$$
See also \cite{OtV, St1, StR, V1, V2} for the extension of Otto's work to Riemannian manifolds.

The Witten Laplacian is a natural extension of the standard Laplace Beltrami operator and has a close connection with differential geometry, probability theory, quantum field theory and statistical mechanics. In view of this, it is natural to raise the question whether one can extend the results which hold for the standard Laplace Beltrami operator to the Witten Laplacian on manifolds. The main tool which leads such an extension possible is the so-called Bakry-Emery Ricci curvature associated to $L$  was introduced in \cite{BE}, i.e.,
\begin{eqnarray*}
Ric(L)=Ric+\nabla^2\phi,
\end{eqnarray*}
which plays the same r\^ole as the Ricci curvature for the standard Laplace Beltrami operator. We refer the reader to \cite{BQ1, BL, Li05} for the Li-Yau Harnack estimates and the heat kernel estimates to the heat equation $(\ref{L-HE})$, and to \cite{Li05} for the  extension of S.-T. Yau's Strong Liouville theorem for the positive $L$-harmonic functions and the $L^1$-uniqueness of the heat equation $(\ref{L-HE})$ on complete Riemannian manifolds. See also \cite{AN, BQ2, FLZ, FLL, OtV, StR, V1, V2, WW} for other results on the study of the Witten Laplacian and Bakry-Emery Ricci curvature on manifolds with weighted measures.

The Bakry-Emery Ricci curvature has been essentially used in Perelman's work on the entropy formula for Ricci flow. In \cite{P1}, Perelman first introduced the $\mathcal{F}$-functional on the space of Riemannian metrics and smooth functions, i.e., $\mathcal{M}=\{ {\rm all\ Riemannian\ metrics}\ g \ {\rm on}\ M\}\times C^\infty(M)$, as follows,
\begin{eqnarray*}
\mathcal{F}(g, f)=\int_M (R+|\nabla f|^2)e^{-f}dv,
\end{eqnarray*}
where $f\in C^\infty(M)$, $R$ denotes the scalar curvature on $(M, g)$, and $dv$ denotes the volume measure.  Under the constraint condition which requires that
\begin{eqnarray*}
d\mu=e^{-f}dv\end{eqnarray*}
is a fixed weighted measure on $(M, g)$, Perelman proved that the gradient flow of $\mathcal{F}$ with respect to the standard $L^2$-metrics on  $\mathcal{M}$ is given by the following modified Ricci flow
\begin{eqnarray*}
\partial_t g=-2(Ric+\nabla^2 f),
\end{eqnarray*}
and $f$ satisfies the so-called conjugate heat equation
\begin{eqnarray*}
\partial_t f=-\Delta f-R.
\end{eqnarray*}
Moreover, Perelman \cite{P1} introduced the $W$-entropy and proved its monotonicity for the Ricci flow on compact manifolds.
This result plays an important r\^ole in the proof of the no local collapsing theorem and in the final resolution of the Poincar\'e conjecture and geometrization conjecture (see e.g. \cite{CZ, MT, KL}). Since then, many people have derived the $W$-entropy formula for various
geometric evolution equations and used it to study further  analysis and geometric properties of manifolds. See e.g. \cite{CLN, CKL, KL, N1,N2, Ec, LNVV, KN}.

In \cite{Li11a}, see also \cite{Li11b, Li12}, inspired by Perelman's work on the $W$-entropy formula for Ricci flow, the second author proved the $W$-entropy formula for the fundamental solution of the Witten-Laplacian on complete Riemannian manifolds with bounded geometry condition, which extends a previous result due to Ni \cite{N1, N2}, who proved an analogue of Perelman's $W$-entropy formula for the heat equation $\partial_t u=\Delta u$ on complete Riemannian manifolds with a fixed metric. More precisely, we have

\begin{theorem}\label{Th-A} (\cite{Li11a, Li11b, Li12})  Let $(M, g)$ be a complete Riemannian manifold with bounded geometry condition\footnote{Here we say that $(M, g)$ satisfies the bounded geometry condition if the Riemannian curvature tensor ${\rm Riem}$ and its covariant derivatives $\nabla^k {\rm Riem}$ are uniformly bounded on $M$, $k=1, 2, 3$.}, and $\phi\in C^4(M)$ with $\nabla\phi\in C_b^3(M)$. Let $m>n$, and  $u={e^{-f}\over (4\pi t)^{m/2}}$ be the fundamental associated with the Witten Laplacian, i.e., the heat kernel to
the heat equation $\partial_t u=Lu$.
Let
\begin{eqnarray*}
H_m(u, t)=-\int_M u\log u d\mu-{m\over 2}(1+\log(4\pi t)).
\end{eqnarray*}
Define
\begin{eqnarray*}
W_m(u, t)={d\over dt}(tH_m(u)).
\end{eqnarray*}
Then
\begin{eqnarray}
W_m(u, t)=\int_M\left[t|\nabla \log u|^2+f-m\right]{e^{-f}\over (4\pi t)^{m/2}}d\mu,\label{W-0}
\end{eqnarray}
and
\begin{eqnarray}
{d W_m(u, t)\over dt}&=&-2\int_M
t \left(\left|\nabla^2 f-{g\over 2t}\right|^2+Ric_{m,
n}(L)(\nabla f,
\nabla f)\right)ud\mu\nonumber\\
& & \hskip2cm  -{2\over m-n}\int_M t \left({\nabla\phi\cdot\nabla
f}+{m-n\over 2t}\right)^2ud\mu,\label{W-1}
\end{eqnarray}
where
\begin{eqnarray*}
Ric_{m, n}(L)=Ric+\nabla^2\phi-{\nabla\phi\otimes\nabla\phi\over m-n}.
\end{eqnarray*}
is the $m$-dimensional Bakry-Emery Ricci curvature associated with the Witten Laplacian $L$.

In particular, if $(M, g, \phi)$ satisfies the bounded geometry condition and $Ric_{m, n}(L)\geq 0$, then the $W$-entropy is decreasing in time $t$, i.e.,
\begin{eqnarray*}
{dW_m(u, t)\over dt}\leq 0, \ \ \ \ \ \ \ \forall t\geq 0.
\end{eqnarray*}
\end{theorem}

The purpose of this paper is to extend the $W$-entropy formula in Theorem \ref{Th-A} to the heat equation $(\ref{L-HE})$ associated with the time dependent Witten Laplacian on compact Riemannian manifolds equipped with time dependent metrics and potentials. In view of Perelman's work using the $W$-entropy formula for the Ricci flow to remove ``the major stumbling block in Hamilton's approach to geometrization" (\cite{P1}),  it might be possible that the $W$-entropy formula for the time dependent Witten Laplacian can  bring some new insights to the study of  geometric analysis on Riemannian manifolds with time dependent metrics and potentials.

We are now in a position to state the main results of this paper as follows.

\begin{theorem}\label{Th-B} Let $(M, g(t), t\in [0, T])$ be a family of compact Riemannian manifolds with potential
functions
$f(t)\in C^\infty(M)$, $t\in [0, T]$. Suppose that $g(t)$ and $f(t)$
satisfy the conjugate
equation
\begin{eqnarray*}
\frac{\partial f} {\partial t}={1\over 2}{\rm Tr}\left(
\frac{\partial g}{\partial t}\right).
\end{eqnarray*}
Let
\begin{eqnarray*}
L=\Delta_{g(t)}-\nabla_{g(t)}f(t)\cdot\nabla_{g(t)}
\end{eqnarray*}
be the time dependent Witten Laplacian on $(M, g(t), f(t))$. Let $u$ be a positive solution of the heat equation
\begin{eqnarray*}
\partial_t u = Lu
\end{eqnarray*}
with initial date $u(0)$ satisfying $\int_M
u(0)d\mu(0)=1$.
Let
\begin{eqnarray*}
H_m(u, t)=-\int_M u\log u d\mu-{m\over 2}(1+\log(4\pi t)).
\end{eqnarray*}
Define
\begin{eqnarray*}
W_m(u, t)={d\over dt}(tH_m(u)).
\end{eqnarray*}
Then
\begin{eqnarray*}
W_m(u, t)=\int_M \left[t|\nabla\log u|^2-\log u-m\right]ud\mu,
\end{eqnarray*}
and
\begin{eqnarray}
{d\over dt}W_m(u, t)&=&-2t\int_M \left|\nabla^2\log u+{g\over 2t}\right|^2ud\mu-{2t\over m-n}\int_M \left(\nabla \phi\cdot \nabla\log u-{m-n\over 2t}\right)^2  ud\mu\nonumber\\
& &\ \ \ \ \ \ \ \ \ -2\int_M t\left({1\over 2}{\partial g\over \partial t}+Ric_{m, n}(L)\right)(\nabla \log u, \nabla \log u)ud\mu.\label{New-W-1}
\end{eqnarray}
In particular, if $\{g(t), f(t), t\in (0, T]\}$ satisfies the $m$-dimensional Perelman's super Ricci flow and the conjugate equation
\begin{eqnarray}
& &{1\over 2}{\partial g\over \partial t}+Ric_{m, n}(L)\geq 0,\label{Cond-1}\\
& &\hskip0.5cm \frac{\partial f} {\partial t}={1\over 2}{\rm Tr}\left(
\frac{\partial g}{\partial t}\right),\label{Cond-2}
\end{eqnarray}
then $W_m(u, t)$ is decreasing in $t\in (0, T]$, i.e.,
\begin{eqnarray*}
{dW_m(u, t)\over dt}\leq 0, \ \ \ \forall t\in (0, T].
\end{eqnarray*}
\end{theorem}

As an application of the $W$-entropy formula for the Witten Laplacian on manifolds with time dependent metrics and potentials,  we prove that the optimal logarithmic Sobolev constant associated with the  Witten Laplacian on compact manifolds equipped with the $m$-dimensional Perelman's super Ricci flow is  decreasing in time. More precisely, we have

\begin{theorem}\label{Th-E} Let $(M, g(t), f(t), t\in [0, T])$ be as in Theorem \ref{Th-B}. Then there exists a positive and smooth function
$u=e^{-{v\over 2}}$ such that $v$ achieves the optimal logarithmic Sobolev constant $\mu(t)$ defined by
\begin{eqnarray*}
\mu(t):=\inf\limits\left\{W_m(u, t): \
\int_M {e^{-v}\over (4\pi t)^{m/2}}d\mu=1\right\}.
\end{eqnarray*}
Indeed, $u=e^{-{v\over 2}}$ is a solution to the nonliniear PDE
\begin{eqnarray*}
-4t L u-2u\log u -m u=\mu(t) u.
\end{eqnarray*}
Moreover, if $\{g(t), f(t), t\in [0, T]\}$ satisfies the $m$-dimensional Perelman's super Ricci flow $(\ref{Cond-1})$ and the conjugate equation $(\ref{Cond-2})$, then  $\mu(t)$ is decreasing in $t$ on $[0, T]$.
\end{theorem}

\begin{remark}
{\rm We believe that, by the approach used in \cite{Li11a, Li11b, Li12}, it would be possible to further extend the $W$-entropy formula in Theorem \ref{Th-B} to the fundamental solution of the heat equation associated with the Witten Laplacian on complete Riemannian manifolds with time dependent metrics and potentials satisfying the bounded geometry condition. Technically, this requires to prove  some Hamilton type gradient estimates for the logarithmic of the heat kernel of the Witten Laplacian on complete Riemannian manifolds with time dependent metrics and potentials satisfying the uniformly  bounded geometry condition\footnote{Here we say that $(M, g(t), f(t), t\in [0, T])$ satisfies the uniformly bounded geometry condition if there exists  some $N\in \mathbb{N}$ such that for all $\varepsilon\in (0, T)$, the $k$-th order covariant derivatives of the Riemannian curvature tensor ${\rm Riem(g(t))}$, i.e., $\nabla^k {\rm Riem(g(t))}$, as well as the $k$-th order covariant derivatives of $f(t)$, i.e., $\nabla^k f(t)$ are uniformly bounded on $[\varepsilon, T]\times M$, $k=0, ..., N$.}. We will study this problem in a forthcoming paper. If this can be verified,  we can derive that, if $\{g(t), f(t), t\in (0, T]\}$ is a family of metrics and potentials satisfying $(\ref{Cond-1})$ and $(\ref{Cond-2})$ on a complete Riemanian manifold $M$ with uniformly bounded geometry condition, then
\begin{eqnarray*}
{d\over dt}W_m(u, t)=0, \ \ \ \ \ \ {\rm for\ \ some}\ \ t=\tau \in (0, T],
\end{eqnarray*}
if and only if
\begin{eqnarray*}
\nabla^2\log u&=&-{g\over 2\tau},\\
\nabla f\cdot \nabla \log u&=&{m-n\over 2\tau}.
\end{eqnarray*}
By the same argument as used in \cite{Li11a, Li11b, Li12}, we can further prove the following rigidity theorem: Let $\{g(t), f(t), t\in (0, T]\}$ be a family of metrics and potentials satisfying $(\ref{Cond-1})$ and $(\ref{Cond-2})$ on a complete Riemanian manifold $M$ with uniformly bounded geometry condition. Let $u$ be the fundamental solution to the heat equation $\partial_t u=Lu$. Then,
 \begin{eqnarray*}
{d\over dt}W_m(u, t)=0, \ \ \ \ \ \ {\rm for\ \ some}\ \ t=\tau\in (0, T],
\end{eqnarray*}
if and only if $(M, g(t))$ is isometric to $\mathbb{R}^n$, $f(t)$ is identically equal to a constant, $m=n$, and
$$u(x, t)={e^{-{\|x\|^2\over 4t}}\over (4\pi t)^{n/2}}, \ \ \ \forall \ x\in M=\mathbb{R}^n, \ \ t\in (0, T].$$
}
\end{remark}

The rest of this paper is organized as follows. In Section $2$, we first give a new proof\footnote{
One of the advantages of our new proof is that it gives a natural geometric interpretation of the third term appeared in the $W$-entropy formula $(\ref{W-1})$. See Remark \ref{rem1}.}  to Theorem \ref{Th-A}. In Section $3$, we prove the dissipation formula of the Bolztmann-Shannon entropy for the heat equation of the Witten Laplacian on compact manifolds with time dependent metrics and potentials.  In Section $4$, we prove Theorem \ref{Th-B} and Theorem \ref{Th-E}. In Section $5$, we use Perelman's $W$-entropy formula for Ricci flow to derive the $W$-entropy formula for the backward  heat equation of the Witten Laplacian on compact Riemannian manifolds equipped with a modified Ricci flow introduced by Lott \cite{Lo2}. In Section $6$, we extend Theorem \ref{Th-A} and Theorem \ref{Th-B} to the case $Ric_{m, n}(L)\geq -K$ and compact $K$-super $m$-dimensional Bakry-Emery-Ricci flow.

\section{A new proof of Theorem \ref{Th-A}}

To give a new proof of Theorem \ref{Th-A}, we first recall some elementary geometric formulas on warped product metrics.

Let $m\in \mathbb{N}$, $m\geq n$. Let
$\widetilde{M}=M\times N$, where $(N, g_N)$ is a compact Riemannian manifold with dimension $q=m-n$. Let $\phi\in C^2(M)$. We consider the following warped product metric on $\widetilde{M}$:
\begin{eqnarray}
\widetilde{g}=g_M\bigoplus e^{-{2\phi\over q}}g_N.\label{WPM}
\end{eqnarray}
Let  $\nu_N$ be the volume measure on $N$. Then the volume measure on $(\widetilde{M}, \widetilde{g})$ is given by
$$dvol_{\widetilde{M}}=e^{-\phi}dvol_M \otimes d\nu_N.$$
Denote
$$d\mu=e^{-\phi}dvol_M.$$
Then
$$dvol_{\widetilde{M}}=d\mu\otimes d\nu_N.$$
Without loss of generality, we may assume that
$$\nu_N(N)=1.$$

Let $\widetilde\Gamma$ be the Christoffel symbol on $(\widetilde{M}, \widetilde{g})$. By direct calculation, we can verify that
\begin{eqnarray*}
\widetilde\Gamma_{ij}^k=\Gamma_{ij}^k,\ \ \
\widetilde\Gamma_{\alpha\beta}^k=q^{-1}g^{kl}\partial_l\phi g_{\alpha\beta}, \ \ \
\widetilde\Gamma_{\alpha\beta}^\gamma=\Gamma_{\alpha\beta}^\gamma,
\end{eqnarray*}
and
\begin{eqnarray*}
\widetilde\Gamma_{ij}^\alpha=0,\ \ \widetilde\Gamma_{i\alpha}^k=0,\ \ \widetilde\Gamma_{i\alpha}^\beta=0.
\end{eqnarray*}
Let $\widetilde\nabla$ be the Levi-Civita connection on  $(\widetilde{M}, \widetilde{g})$. For any $f\in C^2(M)$, using the formula
\begin{eqnarray*}
\widetilde\nabla_{ij}^2 f&=&\partial_i\partial_j f-\widetilde\Gamma_{ij}^k \partial_k f,\\
\widetilde\nabla_{i\alpha}^2 f&=&\partial_i\partial_\alpha f-\widetilde\Gamma_{i\alpha}^k \partial_k f,\\
\widetilde\nabla_{\alpha\beta}^2 f&=&\partial_\alpha\partial_\beta f-\widetilde\Gamma_{\alpha\beta}^k \partial_k f,
\end{eqnarray*}
we have
\begin{eqnarray}
\widetilde\nabla_{ij}^2 f&=&\nabla_{ij}^2 f,\label{Hess-1}\\
\widetilde\nabla_{\alpha\beta}^2 f&=&-q^{-1}g^{kl}\partial_l \phi\partial_k f g_{\alpha\beta},\label{Hess-2}\\
\nabla_{i\alpha}^2 f&=&0\label{Hess-3}.
\end{eqnarray}
Hence
\begin{eqnarray}
\left|\widetilde{\nabla}^2 f-{\widetilde{g}\over 2t}\right|^2&=&
\left|\nabla^2 f-{g\over 2t}\right|^2+\left|\widetilde{\nabla}_{\alpha\beta}^2f-{g_{\alpha\beta}\over 2t}\right|^2\nonumber\\
&=&\left|\nabla^2 f-{g\over 2t}\right|^2+\left|{g^{kl}\partial_l \phi\partial_k f g_{\alpha\beta}\over q}+{g_{\alpha\beta}\over 2t}\right|^2\nonumber\\
&=&\left|\nabla^2 f-{g\over 2t}\right|^2+\left|\left({\nabla\phi\cdot\nabla f\over m-n}+{1\over 2t}\right)g_{\alpha\beta}\right|^2\nonumber\\
&=&\left|\nabla^2 f-{g\over 2t}\right|^2+{1\over m-n}\left({\nabla\phi\cdot\nabla f}+{m-n\over 2t}\right)^2.\label{ccc}
\end{eqnarray}

The following result was obtained in a private discussion by Bing-Long Chen and the second author in January 2006.

\begin{theorem}\label{Th-CL} (\cite{CL}) The Laplace-Beltrami operator on $(\widetilde{M}, \widetilde{g})$ is given by
\begin{eqnarray*}
\Delta_{\widetilde{M}}=L+e^{-{2\phi\over m-n}}\Delta_{N}.
\end{eqnarray*}
\end{theorem}
{\it Proof}. The proof can be given by a  direct calculation. \hfill $\square$

\medskip
\noindent
{\bf A new proof of Theorem \ref{Th-A}}.\ \ To avoid technical issue, we only consider the case of compact manifolds. Let $u={e^{-f}\over (4\pi t)^{m/2}}: M\rightarrow [0, \infty)$ be a positive solution to the heat equation $\partial_t u=Lu$. Then it satisfies the following heat equation on $(\widetilde{M}, \widetilde{g})$
\begin{eqnarray*}
\partial_t u=\Delta_{\widetilde{M}} u.
\end{eqnarray*}
Since $f$ depends only on the variable in the $M$-direction, we have $\widetilde\nabla f=\nabla f$. Therefore the $W$-entropy functional $W_m(u, t)$ defined by $(\ref{W-0})$ coincides with the $W$-entropy functional $\widetilde{W}_{m}(u, t)$ defined on $(\widetilde{M}, \widetilde{g})$ as follows
\begin{eqnarray}
\widetilde{W}_m(u, t)=\int_{\widetilde{M}}\left[t|\widetilde\nabla f|^2+f-m\right]{e^{-f}\over (4\pi t)^{m/2}}dvol_{\widetilde{M}}.\label{WonWP}
\end{eqnarray}
Applying the $W$-entropy formula for the heat equation $\partial_t u=\Delta u$ on compact Riemannian manifolds with fixed metric due to Ni \cite{N1, N2}  to $(\widetilde{M}, \widetilde{g})$,  we have
\begin{eqnarray}
{d\widetilde{W}_m(u, t)\over dt}&=&-2\int_{\widetilde{M}} t \left(\left|\widetilde\nabla^2 f-{\widetilde g\over 2t}\right|^2+\widetilde{Ric}(\widetilde\nabla\log u, \widetilde\nabla\log u)\right)ud\mu dv_N.\label{WP-W-1}
\end{eqnarray}
By $(\ref{ccc})$, we have
\begin{eqnarray}
\left|\widetilde\nabla^2 f-{\widetilde{g}\over 2t}\right|^2=\left|\nabla^2 f-{g\over 2t}\right|^2+{2\over m-n} \left({\nabla\phi\cdot\nabla
f}+{m-n\over 2t}\right)^2.\label{WP-Hess}
\end{eqnarray}
On the other hand, by \cite{Bes, Lo1, Li05}, we have
\begin{eqnarray}
\widetilde{Ric}(\widetilde\nabla\log u, \widetilde\nabla\log u)=Ric_{m, n}(L)(\nabla\log u,
\nabla\log u).\label{Ric-WP}
\end{eqnarray}
From $(\ref{WP-W-1})$, $(\ref{WP-Hess})$ and $(\ref{Ric-WP})$, we obtain $(\ref{W-1})$. This finishes the new proof of Theorem \ref{Th-A}.\hfill $\square$

\begin{remark}\label{rem1} {\rm One of the advantages of the above proof is that: when $m\in \mathbb{N}$ and $m>n$, the quantity ${1\over {m-n}}\left(\nabla \phi\cdot \nabla f+{m-n\over 2t}\right)^2$ appeared in the $W$-entropy formula in Theorem \ref{Th-A} has a natural geometric interpretation. It corresponds to the vertical component of the  quantity $\left|\widetilde\nabla^2 f-{\widetilde{g}\over 2t}\right|^2$  on the warped product manifold $\widetilde{M}=M\times N$ equipped with the metric
\begin{eqnarray*}
\widetilde{g}=g\bigoplus e^{-{2\phi\over m-n}}g_N.\label{WPM1}
\end{eqnarray*}
}
\end{remark}

\begin{remark}{\rm We would like to mention that, after the first version of this paper \cite{LL13} (arxiv arxiv1303.6019) was posted online, N. Charalambous and Z. Lu posted their preprint \cite{CL} in which they used the warped product approach to prove the Li-Yau differential Harnack inequality on complete Riemannian manifolds with weighted volume measure.  We also found a preprint  \cite{GPT} by H. Guo, R. Philipowski and A. Thalmaier, in which they studied the Boltzmann entropy dissipation formula on manifolds with time dependent metrics.}

\end{remark}

\section{Dissipation formula of the Boltzmann-Shannon entropy}

Let $(M, g(t), f(t))$ be as in Theorem \ref{Th-B}. Following \cite{BE, Lo1, Li05}, we introduce the Bakry-Emery Ricci curvature associated with  $L$ as
$$
Ric(L)=Ric+\nabla^2 f.$$

The purpose of this section is to prove the following dissipation formula for the Boltzmann-Shannnon entropy associated with  the Witten Laplacian on manifolds with time dependent metrics and potentials.

\begin{theorem}\label{Th-B1}  Let $u$ be a positive solution to the heat equation $\partial_t u=Lu$.
Let
$$H(u, t)=-\int_M u\log u d\mu$$
be the Boltzmann-Shannon entropy associated with the Witten Laplacian $L$. Then
\begin{eqnarray}
{\partial^2\over \partial t^2} H(u, t)=-2\int_M \left[|\nabla^2\log u|^2+\left({1\over 2}{\partial g\over \partial t}+Ric(L)\right)(\nabla \log u, \nabla \log u)\right]u d\mu.\label{2ndH-3}
\end{eqnarray}

\end{theorem}
{\it Proof}.  By direct calculation, we have
\begin{eqnarray*}
{\partial\over \partial t} H(u, t)=-\int_M \partial_t u(\log u+1)d\mu=-\int_M Lu (\log u+1)d\mu.
\end{eqnarray*}
Integrating by parts yields
\begin{eqnarray*}
{\partial\over \partial t} H(u, t)=\int_M |\nabla\log u|^2_{g(t)} ud\mu,
\end{eqnarray*}which further implies that, as $\partial_t (d\mu)=0$, we have
\begin{eqnarray}
{\partial^2\over \partial t^2}H(u, t)&=&\int_M {\partial\over \partial t}(|\nabla\log u|^2_{g(t)} u)d\mu\nonumber\\
&=&\int_M \left[{\partial \over \partial t}g^{ij}\nabla_i\log u\nabla_j\log u\right]u d\mu+\int_M {\partial \over \partial t}\left[{|\nabla u|^2\over u}\right]_{\rm g(t)\ fixed}d\mu\nonumber\\
&=&\int_M \left[-{\partial \over \partial t}g_{ij}\nabla_i\log u\nabla_j\log u\right]u d\mu+\int_M {\partial \over \partial t}\left[{|\nabla u|^2\over u}\right]_{\rm g(t)\ fixed}d\mu\nonumber\\
&=&\int_M \left(-{\partial g\over \partial t}(\nabla\log u, \nabla u)+{\partial \over \partial t}\left[{|\nabla u|^2\over u}\right]_{\rm g(t)\ fixed}\right)d\mu,\label{aaa}
\end{eqnarray}
where $[\cdot]_{\rm g(t)\ fixed}$ means that the quantity $|\nabla u|^2$ in $[\cdot]$ is defined under a fixed metric $g(t)$.
and we have used the facts  $|\nabla \log u|^2=g^{ij}\nabla_i\log u\nabla_j\log u$ as well as $\partial_t g^{ij}=-\partial_t g_{ij}$.

By the entropy dissipation formula in \cite{BE, Li12}, we have
\begin{eqnarray}
\int_M
{\partial \over \partial t}\left[{|\nabla u|^2\over u}\right]_{\rm g(t)\ fixed}d\mu=-2\int_M \left[|\nabla^2\log u|^2+Ric(L)(\nabla \log u, \nabla \log u)\right]ud\mu.\label{bbb}
\end{eqnarray}Combining $(\ref{aaa})$ and $(\ref{bbb})$, we finish the proof of Theorem \ref{Th-B1}. \hfill $\square$

\medskip

As an easy consequence of Theorem \ref{Th-B}, we have
\begin{corollary}\label{Cor-1} Let $(M, g(t))$ be a closed manifold with a potential $f(t)$. Suppose that $(g(t), f(t))$ satisfies Perelman's super Ricci flow and the conjugate equation, i.e.,
\begin{eqnarray*}
{\partial g\over \partial t}&\geq &-2Ric(L),\\
\frac{\partial f} {\partial t}&=&{1\over 2}{\rm Tr}\left(
\frac{\partial g}{\partial t}\right).
\end{eqnarray*}
Let $u$ be a positive solution to the heat equation $\partial_t u=Lu$.
Then the Boltzmann-Shannon entropy $$H(u, t)=-\int_M u\log ud\mu$$ is concave in time $t$, i.e.,
\begin{eqnarray*}
{d^2\over dt^2}H(u, t)\leq 0.
\end{eqnarray*}
\end{corollary}

\medskip
\section{Proofs of Theorem \ref{Th-B} and Theorem \ref{Th-E}}

Following \cite{Li12}, we introduce
\begin{eqnarray*}
W(u, t)={d\over dt}(tH(u, t)).
\end{eqnarray*}
By direct calculation, we can prove the following

\begin{proposition} We have
\begin{eqnarray*}
W(u, t)=\int_M \left[t|\nabla\log u|^2-\log u\right]ud\mu,
\end{eqnarray*}
and
\begin{eqnarray}
{d\over dt}W(u, t)&=&-2\int_M t\left[|\nabla^2\log u|^2+\left({1\over 2}{\partial g\over \partial t}+Ric(L)\right)(\nabla \log u, \nabla \log u)\right]ud\mu \nonumber\\
& &\hskip3cm +2\int_M |\nabla\log u|^2 ud\mu.\label{www}
\end{eqnarray}
\end{proposition}

\begin{remark}\label{rem2} {\rm From $(\ref{www})$, we can derive that, if
$$
{1\over 2}{\partial g\over \partial t}+Ric(L)-{1\over t}\geq 0.
$$
Then
$$
{d W(u, t)\over dt}\leq 0.
$$
}
\end{remark}

Let
\begin{eqnarray*}
H_m(u, t)=-\int_M u\log u d\mu-{m\over 2}(1+\log(4\pi t)).
\end{eqnarray*}
Following  \cite{P1, N1, Li11a, Li12}, we define $W_m(u, t)$ by the Boltzmann entropy formula
\begin{eqnarray}
W_m(u, t)={d\over dt}(tH_m(u)).\label{BPF}
\end{eqnarray}
We can verify that $W_m(u, t)$ coincides with the expression given in Theorem \ref{Th-B}, i.e.,
\begin{eqnarray*}
W_m(u, t)=\int_M \left[t|\nabla\log u|^2-\log u\right]ud\mu-{m\over 2}(2+\log(4\pi t)).
\end{eqnarray*}
\noindent{\bf Proof of Theorem \ref{Th-B}}. By $(\ref{BPF})$ and $(\ref{2ndH-3})$ in Theorem \ref{Th-B1}, we have
\begin{eqnarray}
{dW_m(u)\over dt}&=&-2\int_M t\left[|\nabla^2\log u|^2+\left({1\over 2}{\partial g\over \partial t}+Ric(L)\right)(\nabla \log u, \nabla \log u)\right]ud\mu\nonumber\\
 & & \ \ \ \ \ \ \ +2\int_M |\nabla\log u|^2 ud\mu-{m\over 2t}. \label{New-W-2}
\end{eqnarray}
Note that
\begin{eqnarray*}
2t|\nabla^2\log u|^2+{m\over 2t}=2t\left|\nabla^2 \log u+{g\over 2t}\right|^2+{m-n\over 2t}-2\Delta \log u.
\end{eqnarray*}
Hence
\begin{eqnarray*}
{d\over dt}W_m(u, t)&=&-{m-n\over 2t}-2t\int_M \left|\nabla^2\log u+{g\over 2t}\right|^2ud\mu+2\int_M |\nabla\log u|^2ud\mu+2\int_M \Delta \log u ud\mu\\
& &-2\int_M t\left({1\over 2}{\partial g\over \partial t}+Ric(L)\right)(\nabla \log u, \nabla \log u)ud\mu.
\end{eqnarray*}
Integrating by part yields
\begin{eqnarray*}
\int_M \Delta \log u ud\mu&=&\int_M (L\log u+\nabla \phi \cdot \nabla \log u) u d\mu\\
&=&-\int_M |\nabla \log u|^2 ud\mu+\int_M \nabla \phi\cdot \nabla \log u u d\mu,
\end{eqnarray*}
whence
\begin{eqnarray*}
{d\over dt}W_m(u, t)&=&-{m-n\over 2t}-2t\int_M \left|\nabla^2\log u+{g\over 2t}\right|^2ud\mu+2\int_M \nabla \phi\cdot \nabla\log u  ud\mu\\
& &-2\int_M t\left({1\over 2}{\partial g\over \partial t}+Ric(L)\right)(\nabla \log u, \nabla \log u)ud\mu.
\end{eqnarray*}
Note that
\begin{eqnarray*}
& &{m-n\over 2t}+2t Ric(L)(\nabla \log u, \nabla \log u)-2 \nabla \phi\cdot \nabla \log u \\
&=&2t Ric_{m, n}(L)(\nabla \log u, \nabla \log u)+{2t\over m-n}\left(\nabla \phi\cdot \nabla \log u-{m-n\over 2t}\right)^2.
\end{eqnarray*}
Hence
\begin{eqnarray*}
{d\over dt}W_m(u, t)&=&-2t\int_M \left|\nabla^2\log u+{g\over 2t}\right|^2ud\mu-{2t\over m-n}\int_M \left(\nabla \phi\cdot \nabla\log u-{m-n\over 2t}\right)^2  ud\mu\\
& &-2\int_M t\left({1\over 2}{\partial g\over \partial t}+Ric_{m, n}(L)\right)(\nabla \log u, \nabla \log u)ud\mu.
\end{eqnarray*}
This proves the $W$-entropy formula in Theorem \ref{Th-B}. The monotonicity result follows.  The  proof of Theorem \ref{Th-B} is completed. \hfill $\square$

\medskip

\noindent{\bf Proof of Theorem \ref{Th-E}}. The proof is similar to
Perelman \cite{P1}. See also  \cite{Li11a}. By definition, we have
\begin{eqnarray}
\mu(t)=\inf\limits_{u}\left\{\int_M \left[4t |\nabla
u|^2-u^2\log u^2-mu^2\right](4\pi
t)^{-{m\over 2}}d\mu\right\},\label{optmu}
\end{eqnarray}
where $\inf\limits$ is taken among all the $u$ such that
$$\int_M (4\pi t)^{-{m\over 2}}u^2d\mu=1.$$
Indeed, $\mu(t)$ is the optimal constant in the following
logarithmic Sobolev inequality: for all $u$ satisfying the above
condition,
\begin{eqnarray*}
\int_M u^2 \log u^2 (4\pi t)^{-{m\over 2}}d\mu\leq \mu(t)+m+4\int_M
t |\nabla u|^2 (4\pi t)^{-{m\over 2}}d\mu.
\end{eqnarray*}
By a similar argument as used in Perelman \cite{P1} and \cite{CZ,
KL, MT}, we can prove that the minimization problem $(\ref{optmu})$
has a non-negative minimizer $u\in H^1(M, \mu)$, which satisfies the
Euler-Lagrange equation
\begin{eqnarray*}
-4t L u-2u\log u -m u=\mu(t) u.
\end{eqnarray*}
By the regularity theory of elliptic PDEs, we have $u\in C^{1,
\alpha}(M)$. By an argument due to Rothaus \cite{Ro}, we can further
prove that $u$ is strictly positive and smooth. Hence $v=-2\log u$
is also smooth. Moreover, as a consequence of Theorem \ref{Th-B}, we
can derive that $\mu(t)$ is a decreasing  function in $t$
on $[0, T]$ provided that $\{g(t), f(t), t\in [0, T]\}$ satisfies the $m$-dimensional Perelman's super Ricci flow $(\ref{Cond-1})$ and the conjugate equation $(\ref{Cond-2})$. The proof of Theorem \ref{Th-E} is completed.  \hfill
$\square$

\begin{remark}\label{rem3} {\rm  Let $m\in \mathbb{N}$ and $m>n$. Let $(N, g_N)$ be a compact Riemniann manifold of dimension $q=m-n$. Let $\mathcal{M}=M\times N$ be the product manifold equipped with the time dependent warped product metric
\begin{eqnarray*}
\widetilde{g}(t)=g(t)\bigoplus e^{-{2f(t)\over m-n}}g_N.\label{WPM2}
\end{eqnarray*}
Similarly to Remark \ref{rem1}, the quantity ${1\over {m-n}}\left(\nabla f\cdot \nabla\log u-{m-n\over 2t}\right)^2$  appeared in the $W$-entropy formula in Theorem \ref{Th-B} has a natural geometric interpretation. It corresponds to the vertical component of the  quantity $\left|\widetilde\nabla^2 \log u+{\widetilde{g}\over 2t}\right|^2$  on $(\mathcal{M}, \widetilde{g}(t))$.
}
\end{remark}

\begin{remark}\label{rem4}{\rm In \cite{P1}, Perelman gave an interpretation of the $W$-entropy using the Boltzmann entropy formula in statistical mechanics. In \cite{Li11a, Li11b}, the second author gave a probabilistic interpretation of the $W$-entropy for the Ricci flow, the heat equation of the Witten Laplacian and for the Fokker-Planck heat equation. Note that, as in \cite{Li11a, Li11b, Li12}, we have
\begin{eqnarray*}
H_m(u, t)=H(u, t)-H(\gamma, t)
\end{eqnarray*}
where $H(u, t)$ is  the Boltzmann-Shannon entropy associated with the heat equation to the Witten Laplacian on $(M, g(t), f(t))$, and $H(\gamma, t)$ is the Boltzmann-Shannon entropy of the Gaussian heat kernel $\gamma(x, t)$ on $\mathbb{R}^m$ if $m\in \mathbb{N}$ with $m\geq n$, i.e.,
\begin{eqnarray*}
\gamma(x, t)={1\over (4\pi t)^{m/2}}e^{-{|x|^2\over 4t}},
 \ \ \ x\in \mathbb{R}^m, t>0.
\end{eqnarray*}
Thus, in view of the definition formula $(\ref{BPF})$, the $W$-entropy $W_m(u, t)$ can be regarded as the by-product of the Boltzmann-Shannon entropy. This gives a probabilistic interpretation of the $W$-entropy $W_m(u, t)$.

On the other hand, similarly to Perelman \cite{P1}, we can also give a heuristical interpretation of the $W$-entropy using the Boltzmann entropy formula in statistical mechanics:  Suppose that there exists a canonical ensemble with a ``density of state
measure'' $g(E)dE$ such that the partition function
$Z_\beta=\int_{\mathbb{R}^+} e^{-\beta E}g(E)dE$ is given by
\begin{eqnarray}
\log Z_\beta=H_m(u, t),
\end{eqnarray}
where $t=\beta^{-1}$. Here, as in \cite{P1}, we do not discuss the issue whether such a ``density of state measure" exists or not. Then, formally applying the Boltzmann entropy formula in statistical mechanics, the thermodynamical entropy of this canonical ensemble is given by
$$S=\log Z_\beta-\beta {\partial \over \partial \beta}\log Z_\beta.$$
Using the fact  ${\partial \over \partial \beta}={\partial\over \partial t}{\partial t\over \partial \beta}=-{1\over \beta^2}{\partial\over \partial t}=-t^2{\partial\over \partial t}$, we can prove
$$S=W_m(u, t).$$
Moreover, formally using the formula in statistical mechanics
\begin{eqnarray*}
{d S\over d\beta}=-\beta {\partial^2\over \partial\beta^2}\log Z_\beta,
\end{eqnarray*}
we can reprove the $W$-entropy formula in Theorem \ref{Th-B}.
}

\end{remark}

\section{The $W$-entropy for the Ricci flow on warped product manifolds}

Let $m\in \mathbb{N}$ and $m\geq n$. Let $\mathbb{T}^q$ be the $q$-dimensional torus with a fixed flat metric given in local coordinates by $\sum^{q}_{i = 1}dx^2_{i}$, where $q=m-n$. Let $\widetilde{M} = M \times \mathbb{T}^q$ be equipped with a time dependent warped product metric
\begin{eqnarray*}
\widetilde{g}(t)=\sum^{n}_{i, j=1}g_{ij}(t)dx^{i}dx^{j} + u(t)^{\frac{2}{q}}\sum^{q}_{\alpha = 1}dx^2_{\alpha}.
\end{eqnarray*}

In \cite{Lo2}, Lott studied the Ricci flow $\widetilde{g}(t)$ on the warped product manifold $\widetilde{M} = M \times \mathbb{T}^q$, which consists of a modified Ricci flow for the Riemannian metric $g(t)$ and a forward heat equation for a potential function $\psi(t)=-\log u(t)$ on the manifold $M$. In this section we use the Perelman's $W$-entropy formula for the Ricci flow $\widetilde{g}(t)$ on the warped product manifold $\widetilde{M}$ to derive the $W$-entropy formula for the backward heat equation associated with the Witten Laplacian $L=\Delta_{g(t)}-\nabla_{g(t)}\psi(t)\cdot\nabla_{g(t)}$ on the compact manifold $M$ equipped with Lott's modified Ricci flow $g(t)$ and the time dependent potential $\psi(t)$.

We first recall Lott's Ricci flow on $\widetilde{M} = M \times \mathbb{T}^q$.  Let $u = e^{-\psi}$. Let $\widetilde{Ric}$ be the Ricci curvature on $(\widetilde{M}, \widetilde{g})$,  and $Ric$ the Ricci curvature on $(M, g)$. By calculation on warped product manifolds \cite{Bes, Lo1, Lo2}, we have
\begin{eqnarray}
\widetilde{Ric}=Ric^{q}_{\psi} + \frac{1}{q}(\Delta \psi - |\nabla \psi|^2)u^{\frac{2}{q}}\sum^{q}_{i = 1}dx^2_{i},\label{PPP1}
\end{eqnarray}
where $Ric^{q}_{\psi}$ is the $m$-dimensional Bakry-Emery Ricci curvature on $(M, g)$ with respect to the potential function $\psi$, i.e.,
$$
Ric^q_{\psi} :=  Ric + Hess \psi - \frac{1}{q}\nabla \psi \otimes \nabla \psi.
$$
See \cite{BE, Li11a, Li11b, Li12}. Below we will also use the notation $Ric_q$ to denote $\widetilde{Ric}$. By $(\ref{PPP1})$, the scalar curvature on $(\widetilde{M}, \widetilde{g})$, denoted by  $R_{q}$,  is given by
\begin{eqnarray*}
R_{q} = R + 2\Delta \psi - \left(1 + \frac{1}{q}\right)|\nabla \psi|^2.
\end{eqnarray*}
By definition, the Ricci flow on $\widetilde{M}$ is defined by
\begin{eqnarray}
\partial_t\widetilde{g}=-2\widetilde{Ric}.
\label{WRF-1}
\end{eqnarray}
According to Lott \cite{Lo2},  the Ricci flow equation $(\ref{WRF-1})$ is equivalent to the following equations
\begin{eqnarray}
\partial_t g &=&-2Ric_{\psi}^q,\label{WRF-2}\\
\partial_t \psi &=& \Delta \psi - |\nabla \psi|^2.\label{WRF-3}
\end{eqnarray}
Note that, the first equation $(\ref{WRF-2})$ is indeed a modified Ricci flow equation for the metric $g(t)$ on $M$, and the second one $(\ref{WRF-3})$ is a forward heat equation for the potential function $\psi(t)$ on $(M, g(t))$. The systems $(\ref{WRF-2})$ and $(\ref{WRF-3})$ are different from Perelman's (modified)  Ricci flow and the conjugate  heat equation  introduced in \cite{P1}, i.e.,
\begin{eqnarray*}
{\partial g\over \partial t}&=&-2(Ric+\nabla^2 f),\\
\frac{\partial f} {\partial t}&=&-\Delta f-R,
\end{eqnarray*}
and are also different from the following $m$-dimensional Perelman's  Ricci flow and the conjugate heat equation
\begin{eqnarray*}
{\partial g\over \partial t}&=&-2(Ric+\nabla^2 f-{\nabla f\otimes \nabla f\over m-n}),\\
\frac{\partial f} {\partial t}&=&-\Delta f+{|\nabla f|^2\over m-n}-R.
\end{eqnarray*}

Let $\phi$
be a positive solution to the conjugate heat equation on $(\widetilde{M},\widetilde{g})$
\begin{eqnarray}
\partial_t \phi= -\Delta_{\widetilde{M}}\phi + R_{q}\phi.\label{PPP4}\end{eqnarray}
Let $\tau\in [0, T]$ be such that
\begin{eqnarray*}
\partial_t \tau=-1.
\end{eqnarray*}
Write
$$\phi= (4\pi \tau)^{-\frac{n+q}{2}}e^{- \eta}.$$
Then
\begin{eqnarray*}
\partial_t \eta &=& -\Delta_{\widetilde{M}}\eta + |\nabla \eta|^{2} - R_{q} + \frac{n+q}{2\tau},
\end{eqnarray*}
Following Perelman \cite{P1},  the $W$-entropy for the Ricci flow $\widetilde{g}(t)$ on the warped product manifold $\widetilde{M}$ is defined as follows
\begin{eqnarray*}
W(\widetilde{g}, \eta, \tau) &=& \int_{\widetilde{M}} [\tau(|\widetilde\nabla \eta|_{\widetilde{M}}^2 + R_q) + \eta -(n+q)]\phi dvol_{\widetilde{M}},
\end{eqnarray*}
where $dvol_{\widetilde{M}}= u dvol_{M}dvol_{\mathbb{T}^q}$ is volume form on $(\widetilde{M}, \widetilde{g})$.

Applying Perelman's $W$-entropy formula for the Ricci flow  \cite{P1} to $(\widetilde{M}, \widetilde{g})$, we have
\begin{eqnarray}
\frac{d}{d\tau}W(\widetilde{g}, \eta, \tau) = - 2\tau \int_{\widetilde{M}} \left|\widetilde{Ric} + \widetilde{Hess}\eta - \frac{\widetilde{g}}{2\tau}\right|_{\widetilde{M}}^2 \phi dvol_{\widetilde{M}}.\label{PPP2}
\end{eqnarray}
By Theorem \ref{Th-CL}, the Laplace-Beltrami on $(\widetilde{M}, \widetilde{g})$ is given by
\begin{eqnarray*}
\Delta_{\widetilde{M}}=L + u^{-\frac{2}{q}}\Delta_{\mathbb{T}^q},
\end{eqnarray*}
where
\begin{eqnarray*}
L=\Delta-\nabla \psi \cdot \nabla.
\end{eqnarray*}
Here $\Delta$ and $\nabla$ are respectively the Laplace-Beltrami  operator and the gradient operator on $(M, g)$.
In the case $\phi$ is a function depending only on the variable of horizontal direction, the conjugate heat equation $(\ref{PPP4})$ turns out to be
   the following backward heat equation associated with the Witten Laplacian on $(M, g(t))$
\begin{eqnarray}
\partial_t \phi = -L\phi + R_{q}\phi.\label{PPP5}
\end{eqnarray}
In this case, $\eta$ is a function depending only on variable in $M$. Thus
\begin{eqnarray*}
W(\widetilde{g}, \eta, \tau) &=& \int_{M\times \mathbb{T}^q} [\tau(|\nabla \eta|^2 + R_q) + \eta - (n+q)]\phi udvol_{M}dvol_{\mathbb{T}^q}\\
&=& \int_{M} \left[\tau(|\nabla \eta|^2 + R + 2\Delta \psi - \left(1 + \frac{1}{q}\right)|\nabla \psi|^2) + \eta - (n+q)\right]\phi d\mu.
\end{eqnarray*}
Here $d\mu=udvol_{M}$, and we assume $vol(\mathbb{T}^q)=1$.
Note that, for any vector field $v$ on $\widetilde{M}$, by $(\ref{Hess-1})$, $(\ref{Hess-2})$ and $(\ref{Hess-3})$, we have
\begin{eqnarray}
\widetilde{Hess}v = Hess v - \frac{1}{q}u^{\frac{2}{q}}\langle \nabla \psi, \nabla v \rangle \sum^{q}_{\alpha = 1}dx^2_{\alpha}.\label{PPP3}
\end{eqnarray}
Substituting $(\ref{PPP1})$ and $(\ref{PPP3})$ into $(\ref{PPP2})$, we have
\begin{eqnarray*}
\frac{d}{d\tau}W(\widetilde{g}, \eta, \tau) &=& - 2\tau \int_{M} \left|\widetilde{Ric} + \widetilde{Hess}\eta  - \frac{\widetilde{g}}{2\tau}\right|_{\widetilde{M}}^2 \phi d\mu\\
&=& - 2\tau \int_{M} \left|Ric_{\psi}^{q} + Hess \eta - \frac{g}{2\tau} + \frac{1}{q}(\Delta \psi - |\nabla \psi|^2 - \langle \nabla \psi, \nabla \eta \rangle - {q\over 2\tau})u^{\frac{2}{q}}\sum^{q}_{\alpha = 1}dx^2_{\alpha}\right|_{\widetilde{M}}^2 \phi d\mu\\
&=& - 2\tau \int_{M} \left[\left|Ric_{\psi}^{q} + Hess \eta - \frac{g}{2\tau}\right|^2 + \frac{1}{q}(\Delta \psi - |\nabla \psi|^2 - \langle \nabla \psi, \nabla \eta \rangle - {q\over 2\tau})^2\right] \phi d\mu.
\end{eqnarray*}

Thus we have proved the following $W$-entropy formula for the backward heat equation associated with the Witten Laplacian on compact manifolds equipped with Lott's modified Ricci flow and time dependent potentials.

\begin{theorem}\label{mono}
Let $(M, g(t), \psi(t))$ be a compact manifold with a family of Remannian metrics $g(t)$ and potentials $\psi(t)$ which satisfy
\begin{eqnarray*}
\partial_t g &=& -2(Ric + Hess \psi - \frac{1}{q}\nabla \psi \otimes \nabla \psi),\\
\partial_t \psi &=& \Delta \psi - |\nabla \psi|^2 .
\end{eqnarray*}
Let $d\mu= e^{-\psi}dvol_{M}$, and $L=\Delta-\nabla\psi\cdot\nabla$.  Let $\phi$ be a positive solution to the backward  heat equation of the Witten Laplacian on $M$, i.e.,
\begin{eqnarray*}
\partial_t \phi = -L\phi + R_{q}\phi.
\end{eqnarray*}
where $R_{q} = R + 2\Delta \psi - (1 + \frac{1}{q})|\nabla \psi|^2$.
Define the W-entropy $W_{q}(g, \psi, \eta, \tau)$ by
\begin{eqnarray*}
W_{q}(g, \psi, \eta, \tau) = \int_{M} [\tau(|\nabla \eta|^2 + R_{q}) + \eta - (n+q)]\phi d\mu.
\end{eqnarray*}
where $\phi=(4\pi \tau)^{-\frac{n+q}{2}}e^{- \eta}$, and $(\eta, \tau)$ satisfies
\begin{eqnarray*}
\partial_t \eta &=& -L\eta + |\nabla \eta|^{2} - R_{q} + \frac{n+q}{2\tau},\\
\partial_t \tau &=& -1.
\end{eqnarray*}
Then
\begin{eqnarray*}
\frac{d}{d\tau}W_{q}(g, \psi, \eta, \tau) = - 2\tau \int_{M} [|Ric_{\psi}^{q} + Hess \eta - \frac{g}{2\tau}|^2 + \frac{1}{q}(\Delta \psi - |\nabla \psi|^2 - \langle \nabla \psi, \nabla \eta \rangle -{q\over 2\tau})^2] \phi d\mu.
\end{eqnarray*}
In particular, $W_{q}(g, \psi, \eta, \tau)$ is decreasing in the backward time $\tau$, and the monotonicity is strict unless that
\begin{eqnarray*}
Ric_{\psi}^{q} + Hess \eta&=&\frac{g}{2\tau},\\
\Delta \psi-|\nabla \psi|^2&=&\langle \nabla \psi, \nabla \eta \rangle-{q\over 2\tau}.
\end{eqnarray*}

\end{theorem}

As an application of the $W$-entropy formula in Theorem \ref{Th-E}, we have the following result.

\begin{theorem}\label{Th-F}
Let $(M, g(t), \psi(t))$ be a compact manifold with a family of Remannian metrics $g(t)$ and potentials $\psi(t)$ which satisfy
\begin{eqnarray*}
\partial_t g &=& -2(Ric + Hess \psi - \frac{1}{q}\nabla \psi \otimes \nabla \psi),\\
\partial_t \psi &=& \Delta \psi - |\nabla \psi|^2 .
\end{eqnarray*}Then there exists a positive and smooth function
$u=e^{-{\eta\over 2}}$ such that $\eta$ achieves the optimal logarithmic Sobolev constant $\mu(\tau)$ defined by
\begin{eqnarray*}
\mu(\tau):=\inf\limits\left\{W_q(g, \psi, \eta, \tau) : \
\int_M {e^{-\eta}\over (4\pi \tau)^{n+q\over 2}}d\mu=1\right\}.
\end{eqnarray*}
where
\begin{eqnarray*}
W_q(g, \psi, \eta, \tau)=
\int_{M} [\tau(|\nabla \eta|^2 + R_{q}) + \eta - (n+q)]\phi d\mu,
\end{eqnarray*}
Indeed, $u=e^{-{\eta\over 2}}$ is a solution to the nonliniear PDE
\begin{eqnarray*}
-4\tau L u+\tau R_qu-2u\log u -(n+q) u=\mu(\tau) u.
\end{eqnarray*}
Moreover,  $\mu(\tau)$ is decreasing in $\tau$ on $[0, T]$.

\end{theorem}
{\it Proof}. The proof is similar to Perelman's monotonicity theorem for the $\mu$-invariant for Ricci flow \cite{P1}.  See also \cite{CZ, CLN, KL, MT} and the proof of Theorem \ref{Th-E}. \hfill
$\square$
\section{The $W$-entropy formula for Witten Laplacian with negative Bakry-Emery Ricci curvature }

The $W$-entropy formula $(\ref{W-1})$ only implies the monotonicity of the $W$-entropy for the Witten Laplacian on complete Riemannian manifolds with  non-negative $m$-dimensional Bakry-Emery Ricci curvature, and the $W$-entropy formula $(\ref{New-W-1})$ only implies the monotonicity of the $W$-entropy for the Witten Laplacian on compact Riemannian manifolds with time dependent metrics and potentials satisfying the super $m$-dimensional Bakry-Emery Ricci flow and the conjugate heat equation. On the other hand, J. Li and X. Xu \cite{LX} introduced a $W$-entropy for the heat equation $\partial_t u=\Delta u$ on complete Riemannian manifolds with Ricci curvature bounded from below by a negative constant. In this section, we combine the ideas in \cite{LX, Li11a, Li12} and Section $4$ to extend Theorem \ref{Th-A} to Witten Laplacian
on complete Riemannian manifolds with $Ric_{m, n}(L)$ bounded from below by a negative constant, and extend Theorem \ref{Th-B} to Witten Laplacian on compact Riemannian manifolds with time dependent metrics and potentials satisfying the $K$-super $m$-dimensional Bakry-Emery Ricci flow and the conjugate heat equation.

Recall the following entropy dissipation formulas on complete Riemannian manifolds. 

\begin{theorem}\label{AAA}(\cite{Li11a, Li12})
Let~$(M,g)$~be a complete Riemannian manifold with bounded geometry condition, and $\phi\in C^4(M)$ with $\nabla\phi\in C_b^3(M)$. Let $u$ be the fundamental solution to the heat equation ~$\partial_tu=L u$.
Let 
$$H(u, t)=-\int_M u\log ud\mu.$$ Then
\begin{eqnarray*}
{d\over dt}H(u, t)=\int_M {|\nabla u|^2\over u}d\mu,
\end{eqnarray*}
and
\begin{eqnarray*}
{d^2\over dt^2}H(u, t)=-2\int_M (|\nabla^2\log u|^2+Ric(L)(\nabla
\log u, \nabla \log u))ud\mu.
\end{eqnarray*} 
\end{theorem}

\begin{proposition}\label{WW1} Let $m\geq n$ and $K\geq 0$ be two constants. Under the same assumption as in Theorem \ref{AAA}, define
\begin{eqnarray*}
H_{m,K}(u, t)=-\int_M u\log u d\mu-{m\over 2}(1+\log(4\pi t))-\frac m2Kt\Big(1+\frac16Kt\Big),
\end{eqnarray*}
Then
\begin{eqnarray*}
{d\over dt}H_{m, K}(u, t)=\int_M \left[{|\nabla u|^2\over u^2}-{m\over 2t}-{mK\over 2}\left(1+{Kt\over 3}\right)\right] u d\mu.
\end{eqnarray*}
In particular, if $Ric_{m, n}(L)\geq -K$, then
\begin{align*}
\frac{d}{dt}H_{m, K}(u, t) \le0,
\end{align*}
\end{proposition}
{\it Proof}. By Theorem \ref{AAA}, we have
\begin{eqnarray*}
{d\over dt}H_{m, K}(u, t)=\int_M \left[{|\nabla u|^2\over u^2}-{m\over 2t}-{mK\over 2}\left(1+{Kt\over 3}\right)\right] u d\mu.
\end{eqnarray*}
 By the same argument as in \cite{LX}, or using the warped product approach as in \cite{CL} and the Li-Yau type differential Harnack inequality obtained by J. Li and X. Xu in \cite{LX}, we can prove that, if $Ric_{m, n}(L)\geq -K$, then
\begin{eqnarray*}
{|\nabla u|^2\over u^2}-\left(1+{2\over 3}Kt\right){\partial_t u\over u}\leq {m\over 2t}+{mK\over 2}\left(1+{Kt\over 3}\right),
\end{eqnarray*}
from which and using the fact $\int_M \partial_t ud\mu=\int_M Lud\mu=0$, we have 
\begin{eqnarray*}
{d\over dt}H_{m, K}(u, t)\leq 0.
\end{eqnarray*}
\hfill $\square$

We now prove the main result of this section. 

\begin{theorem}\label{WW} Let $m\geq n$ and $K\geq 0$ be two constants.  Under the same assumption as in Theorem \ref{AAA}, define
the $W$-entropy by the Boltzmann formula 
\begin{align*}
W_{m, K}(u, t)={d\over dt}(tH_{m,K}(u)).
\end{align*}
Denote $u=\frac{e^{-f}}{(4\pi t)^{m/2}}$. Then
\beq
W_{m, K}(u, t)=\int_M\left(t|\nabla f|^2+f-m\Big(1+\frac12Kt\Big)^2\right)u d\mu,\label{WmK}
\eeq
and
\begin{align}
\frac{d}{dt}W_{m, K}(u, t)
&=-2t\int_M\left(\Big|\nabla^2 f-\left(\frac1{2t}+\frac K2\right) g\Big|^2 +({\rm Ric}_{m,n}(L)+Kg)(\nabla f, \nabla f) \right)u\, d\mu\notag\\
&\hskip2cm -\frac{2t}{m-n}\int_M\left(\nabla \phi\cdot\nabla f+(m-n)\Big(\frac1{2t}+\frac K2\Big)\right)^2u\,d\mu.\label{WmK2}
\end{align}
In particular, if $Ric_{m, n}(L)\geq -K$, then
\begin{align*}
\frac{d}{dt}W_{m, K}(u, t)
\le0.
\end{align*}
\end{theorem}
{\it Proof}. 
By direct calculation, we can prove $(\ref{WmK})$. By Theorem \ref{AAA}, we have
\begin{align*}
{d\over dt}W_{m,K}(u, t)=&{d\over dt}W(u, t)-\frac{m}{2t}-mK\Big(1+\frac K2t\Big)\\
=&-2\int_M t\left[|\nabla^2\log u|^2+Ric(L)(\nabla \log u, \nabla \log u)\right]ud\mu\\
&+2\int_M |\nabla\log u|^2 ud\mu-\frac{m}{2t}-mK\Big(1+\frac K2t\Big).
\end{align*}
Note that
\begin{align*}
&2t|\nabla^2\log u|^2+\frac{m}{2t}+mK\Big(1+\frac K2t\Big)\\
=&2t\Big|\nabla^2\log u+\frac1{2t}g+\frac{K}{2}g\Big|^2-2(1+Kt)\Delta\log u+(m-n)\Big(\frac1{2t}+K\Big(1+\frac K2t\Big)\Big).
\end{align*}
Integrating by part yields 
\begin{align*}
&{d\over dt}W_{m,K}(u, t)\\
=&-2\int_M t\left[\Big|\nabla^2\log u+\frac1{2t}g+\frac{K}{2}g\Big|^2+\left({\rm Ric}(L)+Kg\right)(\nabla \log u, \nabla \log u)\right]ud\mu\\
&+2(1+Kt)\int_M\nabla\log u\cdot\nabla\phi ud\mu+(m-n)\Big(\frac1{2t}+K\Big(1+\frac K2t\Big)\Big)\\
=&-2\int_M t\left[\Big|\nabla^2\log u+\frac1{2t}g+\frac{K}{2}g\Big|^2+\left({\rm Ric}_{m,n}(L)+Kg\right)(\nabla \log u, \nabla \log u)\right]ud\mu\\
&-\frac{2t}{m-n}\int_M\left(\nabla \log u\cdot\nabla \phi-(m-n)\Big(\frac1{2t}+\frac K2\Big)\right)^2u\,d\mu.
\end{align*}
In particular, if ${\rm Ric}_{m,n}(L)\ge-Kg$, $W_{m, K}(u, t)$ is monotone decreasing.
\hfill $\square$

\begin{remark}{\rm Suppose that $Ric_{m, n}(L)\geq -K$. By Theorem \ref{WW}, ${d\over dt}W_{m, K}(u, t)=0$ if and only if
\begin{eqnarray*}
Ric_{m, n}(L)&=&-Kg,\\
\nabla^2 f&=&\left({1\over 2t}+{K\over 2}\right)g,\\
\nabla\phi\cdot\nabla f&=&(m-n)\Big(\frac1{2t}+\frac K2\Big).
\end{eqnarray*}
In particular, if $m=n$ and $\phi=C$ is a constant, then
$(M, g)$ is an Einstein manifold with $Ric=-K$,  and the potential $f$  satisfies
the shrinking gradient Ricci soliton equation (see also \cite{LX})
\begin{eqnarray*}
{1\over 2}Ric+\nabla^2 f={g\over 2t}.
\end{eqnarray*}
}

\end{remark}

\begin{remark}{\rm Similarly to Section $2$, in the case $m\in \mathbb{N}$, $m\geq n$, and $M$ is a compact Riemannian manifold, we can give a new proof of Theorem \ref{WW} via the warped product method. Let $\widetilde{M}=M\times N$, where $(N, g_N)$ is a compact Riemannian manifold with dimension $q=m-n$. Consider the following warped product metric on $\widetilde{M}$:
\begin{eqnarray*}
\widetilde{g}=g_M\bigoplus e^{-{2\phi\over q}}g_N.\label{WPM}
\end{eqnarray*}
Applying the $W$-entropy formula due to J. Li and X. Xu \cite{LX} for the heat equation $\partial_t u=\Delta_{\widetilde{M}} u$ on $(\widetilde{M}, \widetilde{g})$,  we have
\begin{eqnarray}
{d\over dt}\widetilde{W}_{m, K}(u, t)=-2\int_{\widetilde{M}} t \left(\left|\widetilde\nabla^2 f-{\widetilde g\over 2t}-{K\widetilde{g}\over 2}\right|^2+(\widetilde{Ric}+K\widetilde{g})(\widetilde\nabla\log u, \widetilde\nabla\log u)\right)ud\mu dv_N.\label{WPW-1}
\end{eqnarray}
By $(\ref{ccc})$, we have
\begin{eqnarray}
\left|\widetilde\nabla^2 f-{\widetilde{g}\over 2t}-{K\widetilde{g}\over 2}\right|^2=\left|\nabla^2 f-{g\over 2t}-{Kg\over 2}\right|^2+{2\over m-n} \left({\nabla\phi\cdot\nabla
f}+(m-n)\Big(\frac1{2t}+\frac K2\Big)\right)^2.\label{WPHess}
\end{eqnarray}
On the other hand, by \cite{Bes, Lo1, Li05}, we have
\begin{eqnarray}
(\widetilde{Ric}+K\widetilde{g})(\widetilde\nabla\log u, \widetilde\nabla\log u)=(Ric_{m, n}(L)+Kg)(\nabla\log u,
\nabla\log u).\label{Ric-WP}
\end{eqnarray}
From $(\ref{WPW-1})$, $(\ref{WPHess})$ and $(\ref{Ric-WP})$, we reprove $(\ref{WmK2})$. Note that $(\ref{WPHess})$  also gives us a natural geometric interpretation of the third term in the $W$-entropy formula $(\ref{WmK2})$.
}
\end{remark}

We now extend Theorem \ref{WW} to the Witten Laplacian on compact manifolds with time dependent metrics and potentials. 

\bth\label{Th-BK}  Let $m\geq n$ and $K\geq 0$ be two constants. Under the same assumption as in Theorem \ref{Th-B}, define
\begin{eqnarray*}
H_{m,K}(u, t)=-\int_M u\log u d\mu-{m\over 2}(1+\log(4\pi t))-\frac m2Kt\Big(1+\frac16Kt\Big),
\end{eqnarray*}
and
\begin{align*}
W_{m,K}(u, t)={d\over dt}(tH_{m,K}(u)).
\end{align*}
Denote $u=\frac{e^{-f}}{(4\pi t)^{m/2}}$. Then
\begin{eqnarray*}
W_{m,K}(u, t)=\int_M\left[t|\nabla f |^2+f-m\Big(1+\frac12Kt\Big)^2\right]{e^{-f}\over (4\pi t)^{m/2}}d\mu,\label{W-0}
\end{eqnarray*}
and
\begin{align*}\label{WMK}
{d\over dt}W_{m,K}(u, t)=&-2 t\int_M\Big|\nabla^2 f-\frac1{2t}g-\frac{K}{2}g\Big|^2{e^{-f}\over (4\pi t)^{m/2}}d\mu\notag\\
&\ \ -2 t\int_M\left({1\over 2}{\partial g\over \partial t}+{\rm Ric}_{m,n}(L)+Kg\right)(\nabla f, \nabla f){e^{-f}\over (4\pi t)^{m/2}}d\mu\notag\\
&\ \ \ \ -\frac{2t}{m-n}\int_M\left(\nabla \phi\cdot\nabla f+(m-n)\Big(\frac1{2t}+\frac K2\Big)\right)^2{e^{-f}\over (4\pi t)^{m/2}}d\mu.
\end{align*}
In particular, if $\{g(t), \phi(t), t\in (0, T]\}$ is the $K$-super $m$-dimensional Bakry-Emery Ricci flow and satisfies the conjugate equation
\begin{eqnarray}
& &{1\over 2}{\partial g\over \partial t}+{\rm Ric}_{m, n}(L)\geq -Kg,\label{Cond-1K}\\
& &\hskip0.5cm \frac{\partial \phi} {\partial t}={1\over 2}{\rm Tr}\left(
\frac{\partial g}{\partial t}\right),\label{Cond-2K}
\end{eqnarray}
then $W_{m,K}(u, t)$ is decreasing in $t\in (0, T]$, i.e.,
\begin{eqnarray*}
{d\over dt}W_{m,K}(u, t)\leq 0, \ \ \ \forall t\in (0, T].
\end{eqnarray*}
\eeth
\proof By \eqref{www}, and replacing $Ric(L)$ by ${1\over 2}{\partial g\over \partial t}+Ric(L)$, the proof is similar to the one of Theorem \ref{WW}. \hfill $\square$

As an application of  Theorem \ref{Th-BK}, we have the following

\begin{theorem}\label{Th-EK} Let $(M, g(t), \phi(t), t\in [0, T])$ be as in Theorem \ref{Th-BK}. Then there exists a positive and smooth function
$u=e^{-{v\over 2}}$ such that $v$ achieves the optimal logarithmic Sobolev constant $\mu_K(t)$ defined by
\begin{eqnarray*}
\mu_K(t):=\inf\limits\left\{W_{m,K}(u, t): \
\int_M {e^{-v}\over (4\pi t)^{m/2}}d\mu=1\right\}.
\end{eqnarray*}
Indeed, $u=e^{-{v\over 2}}$ is a solution to the nonliniear PDE
\begin{eqnarray*}
-4t L u-2u\log u -m\Big(1+\frac K2t\Big)^2 u=\mu_K(t) u.
\end{eqnarray*}
Moreover, if $\{g(t), \phi(t), t\in [0, T]\}$ satisfies the $K$-super $m$-dimensional Bakry-Emery Ricci flow $(\ref{Cond-1K})$ and the conjugate equation $(\ref{Cond-2K})$, then  $\mu_K(t)$ is decreasing in $t$ on $[0, T]$.
\end{theorem}
{\it Proof}. The proof is similar to the one of Theorem \ref{Th-E}. \hfill $\square$

\medskip

\noindent{\bf Acknowledgement}. This paper is an extended version of the previous preprint \cite{LL13} in arxiv (arxiv1303.6019). We would like to thank Prof. D. Bakry, Prof. B.-L. Chen  and Dr. Yu-Zhao Wang for useful discussions. We would like to thank Prof. Zhiqin Lu for the citation of the first version of this paper in \cite{CL}.

\medskip

\begin{flushleft}

Songzi Li\\

\medskip

{\sc School of Mathematical Science, Fudan University, \ 220, Handan Road, Shanghai, 200432, China}\\

{\sc Institut de Math\'ematiques de Toulouse, Universit\'e Paul Sabatier, 118, route de Narbonne, 31062, Toulouse Cedex 9, France}

\medskip

Xiang-Dong Li\\

\medskip

{\sc  Academy of Mathematics and Systems Science, Chinese
Academy of Sciences, 55, Zhongguancun East Road, Beijing, 100190, China}, E-mail: xdli@amt.ac.cn
\end{flushleft}

\end{document}